\title{A note on the double-critical graph conjecture}
\author{Hao Huang and Alexander Yu\\
			Department of Mathematics and Computer Science\\
			Emory University\\
			Atlanta, GA 30332}
\date{}

\author{Hao Huang
	\thanks{Department of Mathematics and Computer Science, Emory University, Atlanta, GA 30322. Email: hao.huang@emory.edu.}
	\and Alexander Yu
	\thanks{Department of Mathematics and Computer Science, Emory University, Atlanta, GA 30322. Email: alexyuyuyu@gmail.com.}
}

\documentclass[12]{article}
\usepackage{amssymb}
\usepackage{color}

\newtheorem{theorem}{\bf Theorem}[section]
\newtheorem{lemma}[theorem]{Lemma}

\newtheorem{conjecture}[theorem]{Conjecture}
\newenvironment{pf}{\noindent {\em Proof.}}{\hfill\rule{2mm}{3mm}\par\medskip}

\def\qed{\hfill \rule{4pt}{7pt}}

\def\pf{\noindent {\it Proof. }}

\begin{document}

\maketitle

\begin{abstract}
A connected $n$-chromatic graph $G$ is double-critical if for all the edges $xy$ of $G$, the graph $G-x-y$ is $(n-2)$-chromatic. In 1966, Erd\H os and Lov\'asz conjectured that the only double-critical $n$-chromatic graph is $K_n$. This conjecture remains unresolved for $n \ge 6.$ In this short note, we verify this conjecture for claw-free graphs $G$ of chromatic number $6$.
\end{abstract}

\section{Introduction}

In this note, we consider finite and simple graphs.
For a graph $G$, we use $V(G)$ and $E(G)$ to denote the set of vertices and edges of $G$, respectively. A {\it subgraph} of $G$ is a graph whose vertex set is a subset of $V(G)$ and whose edge set is a subset of $E(G)$. We say that a subgraph $H$ is an {\it induced subgraph} of $G$ if, for any $x,y \in V(H)$, $xy \in E(H) \mbox{ iff } xy \in E(G)$.
Let $G$ be a graph and $S \subset V(G).$ Then $G[S]$, the subgraph of $G$ induced by $S,$ denotes the graph with vertex set $S$ and edge set $\{uv \in E(G): u,v \in S\}$, and let $G - S=G[(V(G)\setminus S].$ When $S= \{x,y\}$, we often write $G-x-y$ instead of $G-S$. For a positive integer $k$, a {\it proper $k$-coloring} of a graph $G$ is a function $c$ from $V(G)$ to a set of $k$ colors such that $c(u) \ne c(v)$ for any $uv \in E(G).$ A graph $G$ is {\it $k$-colorable} if $G$ has a proper $k$-coloring. We use $\chi(G)$ to denote the smallest integer $k$ such that $G$ is $k$-colorable, which is known as the {\it chromatic number} of $G$. Further we denote by $\omega(G)$ and $\alpha(G)$ the size of the largest clique and independent set in $G$, respectively, and $N(v)$ the set of vertices adjacent to the vertex $v$ in $G$.


We say that a connected graph $G$ with chromatic number $n$ is {\it $n$-double-critical}, if, for any $xy \in E(G)$, $\chi(G-x-y) = n-2$. It is easy to see that the complete graph $K_n$ is $n$-double-critical. The following elegant conjecture was posed by Erd\H{o}s and Lov\'asz \cite{Er66} more than fifty years ago.
\begin{conjecture} \label{DCGC}
$K_n$ is the only $n$-double-critical graph. 
\end{conjecture}

It is easy to see that Conjecture \ref{DCGC} holds for $n \le 3.$ With some extra work it can also be verified for $n=4.$ In 1986, Stiebitz \cite{St86} showed that Conjecture \ref{DCGC} is true for $n=5.$ He proved the existence of $K_4$ in every 5-double-critical graph by considering uniquely 3-colorable subgraphs of $G$.
%
%
%
%
%
%
%
%
%
%
%
%
%
%
%
%
%
%
However, this technique does not seem to generalize to finding a larger clique in an $n$-double-critical graph with $n \ge 6.$ 

The double-critical graph conjecture is a special case of a more general conjecture, the so-called Erd\H os-Lov\'asz Tihany conjecture \cite{Er66}: for any graph $G$ with $\chi(G) > \omega(G)$ and any two integers $k, l \ge 2$ with $k+l = \chi(G)+1$, there exists a partition $(S, T)$ of the vertex set such that $\chi(G[S]) \ge k$ and $\chi(G[T]) \ge l$. The Erd\H os-Lov\'asz Tihany conjecture was proved for various cases: $(k, l)=(2,2),(2,3),(2,4),(3,3),(3,4),(3,5)$ (see: \cite{BJ99, Mo86, St86, St87}). Kostochka and Stiebitz \cite{KM08} showed that it is true for all $(k, l)$ for line graphs. Balogh et al. \cite{BKPS09} generalized this result to quasi-line graphs (a graph is a quasi-line graph if the neighbors of every vertex $v$ can be expressed as the union of two cliques) and graphs of independence number $2$. 

A {\it claw} is a $4$-vertex graph with one vertex of degree 3 and the others of degree 1. For convenience, we write $(v;v_1,v_2,v_3)$ to denote a claw in which $v$ has degree 3. A graph $G$ is  {\it claw-free} if it does not have a claw as an induced subgraph. Note that the graphs from both families in Balogh et al.'s result \cite{BKPS09} are claw-free. It would be interesting to know whether the double-critical graph conjecture, or the Erd\H os-Lov\'asz Tihany conjecture holds for all the claw-free graphs. As a step in that direction, we prove the following theorem.


\begin{theorem}\label{thm_main}
Let $G$ be a double-critical graph with $\chi(G)=6$. If $G$ is claw-free, then $G \cong K_6$.
\end{theorem}

The rest of this note is organized as follows. In the next section we will prove several lemmas that will be repeatedly used throughout the proof of Theorem \ref{thm_main}. Section \ref{sec_main} contains the proof of Theorem \ref{thm_main}.

\section{Lemmas}

Given a graph $G$ of chromatic number $n$, and a proper $n$-coloring of $G$, all vertices of the same color form a {\it color class}. By definition, each color class is an independent set in $G$.

\begin{lemma}\label{clique}
Let $G$ be an $n$-double-critical graph. If $\omega(G) \ge n-1$, then $G \cong K_n$. Thus if $G \not\cong K_n$, then $\omega(G) \le n-2$.
\end{lemma}
\pf
Suppose $\omega(G) \ge n-1$. Let $v_1,\ldots, v_{n-1}\in V(G)$ be the vertices that induce a copy of $K_{n-1}$. Among all the proper $n$-colorings of $G$ with color classes $V_1, \ldots, V_n$, and $v_i \in V_i$ for $1 \le i \le n-1$, we choose one that minimizes $|V_n|$.

Let $v_n \in V_n$. We claim that $N(v_n) \cap V_i \ne \emptyset$ for $i=1,\ldots,n-1$.
Suppose not. Without loss of generality we may assume that $N(v_n) \cap V_1 = \emptyset$. If $V_n = \{v_n\}$ then the independent sets $V_1 \cup \{v_n\}, V_2, \ldots ,V_{n-1}$ form an proper $(n-1)$-coloring of $G$, contradicting the assumption that $\chi(G)=n$. If $V_n \setminus \{v_n\} \ne \emptyset$, then $V_1 \cup \{v_n\}, V_2, \ldots ,V_n \setminus \{v_n\}$ are the color classes of an $n$-coloring of $G$. Thus we have a contradiction to the minimality of $|V_n|$.

We now show that $v_n$ is adjacent to every vertex in $\{v_1, \cdots, v_{n-1}\}$. Suppose not. Without loss of generality assume that $v_1 \not\in N(v_n)$.
Then, by the above claim, $v_n$ is adjacent to some $y \in V_1 \setminus \{v_1\}$.
However $\chi(G-v_n-y)=n-2$, a contradiction, since $\{v_1,\ldots,v_{n-1}\}$ induces a copy of $K_{n-1}$ in $G-v_n-y$.

Hence $\{v_1,\ldots,v_n\}$ induces $K_n$ in $G$.
If $V(G) = \{v_1,\ldots,v_n\}$, then $G \cong K_n$.
Suppose $V(G) \ne \{v_1,\ldots,v_n\}$, and let $x \in V(G)$ such that $x \notin \{v_1,\ldots,v_n\}$.
Since $G$ is connected, there exists $z \in V(G)$ such that $xz \in E(G)$. However, $G-x-z$ contains a clique on $n-1$ vertices, which contradicts that $\chi(G-x-z)=n-2$, and thus $G \cong K_n$.
\qed
\bigskip

\begin{lemma}\label{existence}
Let $G$ be an $n$-double-critical graph that is claw-free. For $xy \in E(G)$, let $V_1, \ldots, V_{n-2}$ be the color classes of an $(n-2)$-coloring of $G-x-y$. Then $N(x) \cap N(y) \cap V_i \ne \emptyset$ for $i \in \{1,\ldots,n-2\}$. 
\end{lemma}

\pf
Suppose not. Then, without loss of generality, we may assume that $N(x) \cap N(y) \cap V_{n-2} = \emptyset$. Let $V_x = \{x\} \cup (V_{n-2} \setminus N(x))$ and $V_y = \{y\} \cup V_{n-2} \setminus V_x$. Note that $V_x$ and $V_y$ are independent sets. Now $V_1,\ldots,V_{n-3},V_x,V_y$ are the color classes of an $(n-1)$-coloring of $G$, contradicting the assumption that $\chi(G)=n$.
\qed


\bigskip

The {\it degree} of a vertex $v$, denoted $d(v)$, in a graph is the number of edges incident to $v$. We denote by $\Delta(G)$ and $\delta(G)$ the maximum and minimum degree of a vertex in $G$ respectively. The following lemma shows that for $\chi(G)=6$, it suffices to consider 6-double-critical graphs in which every pair of adjacent vertices has $4$ or $5$ common neighbors.

\begin{lemma}\label{neighbor}
Let $G$ be a $6$-double-critical graph that is also claw-free. If $G \not\cong K_6$, then for any $xy \in E(G), 4 \le |N(x) \cap N(y)| \le 5$. If, in addition, $|N(x) \cap N(y)| =5$, then $G[N(x) \cap N(y)] \cong C_5$.
\end{lemma}

\pf
We may assume that $\Delta(G[N(x)\cap N(y)])\le 2$.
For, otherwise, let $v \in N(x)\cap N(y)$ and let $v_1,v_2,v_3 \in N(x)\cap N(y)$ such that $vv_i \in E(G), i=1,2,3$.
Since $(v;v_1,v_2,v_3)$ does not induce a claw in $G$, there exist $i,j \in \{1, 2, 3\}$ such that $i \ne j$ and $v_iv_j \in E(G)$.
Thus $\{v,v_i,v_j,x,y\}$ induces a copy of $K_5$ in $G$.
Hence by Lemma \ref{clique}, $G \cong K_6$.

We may also assume that $\delta(G[N(x)\cap N(y)]) \ge |N(x)\cap N(y)|-3$.
For, otherwise, let $v \in N(x)\cap N(y)$ and let $v_1,v_2,v_3 \in N(x)\cap N(y)$ such that $vv_i \notin E(G), i=1,2,3$.
Then $v_iv_j \in E(G)$ for all distinct $i,j \in \{1,2,3\}$, since $(v;v_i,v_j,x)$ does not induce a claw in $G$.
Therefore $\{v_1,v_2,v_3,x,y\}$ induces a copy of $K_5$ in $G$.
Once again by Lemma \ref{clique}, $G \cong K_6$.

Hence, if $G \not\cong K_6$, $|N(x)\cap N(y)|-3 \le \delta(G[N(x)\cap N(y)]) \le \Delta(G[N(x)\cap N(y)]) \le 2$.
Thus $|N(x)\cap N(y)| \le 5$. On the other hand, by Lemma \ref{existence}, $|N(x)\cap N(y)| \ge 6-2 = 4$.

Now suppose $|N(x) \cap N(y)| =5$. Then $2=|N(x)\cap N(y)|-3 \le \delta(G[N(x)\cap N(y)]) \le \Delta(G[N(x)\cap N(y)]) \le 2$.
Hence $\delta(G[N(x)\cap N(y)]) = \Delta(G[N(x)\cap N(y)])=2$, so every vertex in $G[N(x)\cap N(y)]$ has degree $2$. The only $2$-regular graph on $5$ vertices is $C_5$, thus $G[N(x)\cap N(y)] \cong C_5$.
\qed

\bigskip
The following lemma is an easy consequence of $G$ being claw-free.

\begin{lemma}\label{easy}
Let $G$ be a claw-free graph, and $S$ an independent set of $G$. Suppose $x \in V(G) \setminus S$. Then $|N(x) \cap S| \le 2.$
\end{lemma}

\pf
Suppose $|N(x) \cap S| \ge 3.$ Let $x_1,x_2,x_3 \in N(x) \cap S.$ Since $S$ is an independent set, $(x;x_1,x_2,x_3)$ induces a claw in $G$, a contradiction.  
\qed


\section{Proof of the Main Result} \label{sec_main}

Theorem \ref{thm_main} follows from the two lemmas in this section. From Lemma \ref{neighbor}, we may assume that the number of common neighbors of any two adjacent vertices is either $4$ or $5$. The first lemma settles the case when there exists a pair with $4$ common neighbors.

\begin{lemma}\label{four}
Let $G$ be a $6$-double-critical graph that is claw-free. If $|N(x) \cap N(y)| =4$ for some $xy \in E(G)$, then $G \cong K_6$.
\end{lemma}

\pf
For an arbitrary $xy \in E(G)$, by Lemma \ref{neighbor}, we have $|N(x)\cap N(y)|\ge4$. Thus $d(x) \ge 5$ and $d(y) \ge 5$. Moreover, if $V_1,V_2,V_3,V_4$ denote the color classes of a 4-coloring of $G-x-y$, it follows from Lemma \ref{easy} that $|N(x) \cap V_i| \le 2$ and $|N(y) \cap V_i| \le 2$ for $i \in \{1,2,3,4\}$. Thus $d(x) \le 9$ and $d(y) \le 9$.\\

\textbf{Claim 1.} If $xy \in E(G)$ and $|N(x) \cap N(y)| = 4$, then $d(x),d(y) \in \{7,8\}$.

Let $N(x)\cap N(y)= \{v_1,v_2,v_3,v_4\}$, and let $V_1,V_2,V_3,V_4$ be the color classes of a 4-coloring of $G-x-y$. By Lemma \ref{existence}, we may assume $v_i \in V_i$, $i=1,2,3,4$. 

Suppose $d(x) \in \{5,6\}.$ Then we may assume that $N(x)\cap V_i = {v_i}$ for $i \in \{2,3,4\}$. Since $|N(x)\cap N(v_1)| \ge 4$, $v_1v_i \in E(G)$ for $i \in \{2,3,4\}$.  If, for every pair of distinct $i,j \in \{2,3,4\}$, $v_iv_j \not\in E(G)$, then $(v_1;v_2,v_3,v_4)$ induces a claw in $G$. Otherwise, there exist distinct $i,j \in \{2,3,4\}$ such that $v_iv_j \in E(G)$. Then $G[\{v_1,v_i,v_j,x,y\}] \cong K_5$. Hence $G \cong K_6$ by Lemma \ref{clique}.

Now suppose that $d(x) = 9.$ 
Let $N(x) \setminus \{v_1,v_2,v_3,v_4,y\} = \{u_1,u_2,u_3,u_4\}$. By Lemma \ref{easy}, we may assume $u_i \in V_i$ for $i\in \{1,2,3,4\}$.
For any distinct $j,k \in \{1,2,3,4\}$, $u_ju_k \in E(G)$; for otherwise $(x;u_j,u_k,y)$ induces a claw. 
Thus $G[\{x,u_1,u_2,u_3,u_4\}] \cong K_5$.
Hence, by Lemma \ref{clique}, $G \cong K_6$.\\

\textbf{Claim 2.} If $xy \in E(G)$ and $|N(x) \cap N(y)| = 4$, then $d(x)=d(y)=8$. 

Let $N(x)\cap N(y)= \{v_1,v_2,v_3,v_4\}$, and let $V_1,V_2,V_3,V_4$ be the color classes of a 4-coloring of $G-x-y$. By Lemma \ref{existence}, we may assume $v_i \in V_i, i=1,2,3,4$. 

Suppose $d(x) = 7$ 
and, by Lemma \ref{easy} and by symmetry, let $u_i \in N(x)\cap V_i \setminus v_i$, $i \in \{1,2\}$.
Since $|N(x)\cap N(u_1)| \ge 4$ and $u_1 \not\in N(y)$, $u_1u_2,u_1v_2,u_1v_3,u_1v_4 \in E(G)$.
Similarly, since $|N(x)\cap N(u_2)| \ge 4$ and $u_2 \not\in N(y)$, $u_2v_1,u_2v_3,u_2v_4 \in E(G)$.

We claim that $y$ does not have a neighbor in $V_1\setminus \{v_1\}$ or $V_2 \setminus \{v_2\}$.
For otherwise, suppose there exists $w_1 \in V_1 \setminus {v_1}$ such that $yw_1 \in E(G)$.
Since $|N(y)\cap N(w_1)| \ge 4$, we have $|N(w_1) \cap \{v_2,v_3,v_4\}| \ge 2$. (This is because $d(y) \le 8$, so $y$ has at most two neighbors not from $\{x, v_1, v_2, v_3, v_4, w_1\}$. If $w_1$ is adjacent to at most one vertex from $\{v_2, v_3, v_4\}$, then $|N(y) \cap N(w_1)| \le 3$.) Similarly since $|N(x)\cap N(v_1)| \ge 4$,  we have $|N(v_1) \cap \{v_2,v_3,v_4\}| \ge 2$. 
Thus there exists $i \in \{2,3,4\}$ such that $v_i \in N(v_1) \cap N(w_1)$. 

Note $(v_i;v_1,u_1,w_1)$ induces a claw in $G$, a contradiction.
Therefore by Claim 1 and Lemma \ref{easy}, $d(y)=7$ and $|N(y) \cap V_i| = 2$ for $i \in \{3,4\}$.
Let $w_i \in N(y) \cap V_i \setminus \{v_i\}$ for $i \in \{3,4\}$.
Since $|N(y)\cap N(w_3)| \ge 4$ and $w_3 \not\in N(x)$, $w_3w_4,w_3v_1,w_3v_2,w_3v_4 \in E(G)$,
and similarly since $|N(y)\cap N(w_4)| \ge 4$ and $w_4 \not\in N(x)$, $w_4v_1, w_4v_2, w_4v_3  \in E(G)$.

We may assume that $v_3v_4 \not\in E(G)$; otherwise $G[\{x,u_1,u_2,v_3,v_4\}] \cong K_5$ and, hence, $G \cong K_6$ by Lemma \ref{clique}. Similarly we may assume $v_1v_2 \not\in E(G)$, otherwise $G[\{y,v_1,v_2,w_3,w_4\}] \cong K_5$ and once again $G \cong K_6$ by Lemma \ref{clique}.
Since $|N(x)\cap N(v_1)| \ge 4$ and $v_1v_2 \not\in E(G)$, $v_1v_3,v_1v_4 \in E(G)$.
Similarly since $|N(x)\cap N(v_2)| \ge 4$ and $v_1v_2 \not\in E(G)$, $v_2v_3,v_2v_4 \in E(G)$. 

We claim that $u_iw_j \in E(G)$ for all $i \in \{1,2\}$ and $j \in \{3,4\}$. Suppose not. By symmetry, we assume $w_4u_1 \not\in E(G)$. 
Since $|N(v_2)\cap N(w_4)| \ge 4$, from the known adjacencies we so far only have $w_3, v_3, y \in N(v_2) \cap N(w_4)$, therefore there exists $w_1 \in V_1 \cup V_3$ such that $w_1v_2,w_1w_4 \in E(G).$
In fact, $w_1 \in V_1$, otherwise then $(w_4;v_3,w_1,w_3)$ induces a claw in $G$, a contradiction. Since $w_4u_1 \not\in E(G)$, $w_1 \ne u_1$.
%
%
%
%
Note that $w_1v_3 \not\in E(G)$ otherwise $(v_3;v_1,u_1,w_1)$ induces a claw and similarly $w_1v_4 \not\in E(G)$ otherwise $(v_4;v_1,u_1,w_1)$ induces a claw.
Then $(v_2;w_1,v_3,v_4)$ induces a claw in $G$, a contradiction since $G$ is claw-free.
Hence, $G[\{u_1,u_2,w_3,w_4\}] \cong K_4$.

Consider the graph $G-x-v_3$. Since $G$ is 6-double-critical, $\chi(G-x-v_3)=4$. Let $c:V(G)\rightarrow \{1,2,3,4\}$ be a 4-coloring of $G-x-v_3$. Since $G[\{u_1,u_2,w_3,w_4\}] \cong K_4$, we may assume $c(u_1)=1$, $c(u_2)=2$, $c(w_3)=3$, $c(u_4)=4$. Then $c(v_1)=1$, since $v_1u_2,v_1w_3,v_1w_4 \in E(G)$. Similarly $c(v_2)=2$ and $c(v_4)=4$. Then, since $yv_1,yv_2,yw_3,yw_4 \in E(G)$, $y$ cannot be colored by any of the colors 1, 2, 3, 4. Thus $G-x-v_3$ is not 4-colorable, a contradiction.
Therefore $d(x) = 8$. Similarly $d(y) =8$.
This completes the proof of Claim 2. \\

Now let us fix $xy \in E(G)$ with $|N(x) \cap N(y)| = 4$. Let $N(x) \cap N(y) = \{v_1,v_2,v_3,v_4\}$ and let $V_1,V_2,V_3,V_4$ be the color classes of a 4-coloring of $G-x-y$. By Lemma \ref{existence}, we may assume $v_i \in V_i$ for $i \in \{1,2,3,4\}$.
By Claim 2, $d(x)=8$. Thus let $w_i \in N(x)\cap V_i \setminus \{v_i\}$, $i \in \{1,2,3\}$.
Note that $G[\{w_1,w_3,w_3\}]$ is a clique of size 3; otherwise suppose, by symmetry $w_1w_2 \not\in E(G)$, then $(x;y,w_1,w_2)$ induces a claw in $G$.\\

\textbf{Claim 3.} For $i \in \{1,2,3\}$, $\{v_1,v_2,v_3,v_4\} \setminus \{v_i\} \subset N(w_i)$.

Suppose otherwise. Without loss of generality, we may assume $\{v_2,v_3,v_4\} \not\subset N(w_1)$. Since $|N(x)\cap N(w_1)| \ge 4$, $|N(w_1)\cap \{v_2,v_3,v_4\}| = 2$, 
and $|N(x) \cap N(w_1)| =4$.
By Claim 2, it suffices to consider the case $d(w_1)=8.$
However, from Lemma \ref{easy}, $$d(w_1) \le 1+|N(w_1)\cap V_2|+|N(w_1)\cap V_3|+|N(w_1)\cap V_4| \le 7,$$
a contradiction. Hence we have Claim 3.\\





Note that for $i \in \{1,2,3\}$, $v_iv_4 \not\in E(G)$ otherwise $\{x,v_i,v_4,w_1,w_2\}$ induces a $K_5$. To avoid the claw $(y; v_i, v_j, v_4)$, we must have $v_iv_j \in E(G)$ for distinct $i, j \in \{1,2,3\}$. In this case $\{x, y, v_1, v_2, v_3\}$ induces a copy of $K_5$, and hence $G \cong K_6$ and the proof of Lemma \ref{four} is complete.
\qed

\medskip
Our next lemma settles the remaining case when every pair of adjacent vertices have exactly $5$ common neighbors.

\begin{lemma}\label{five}
Let $G$ be $6$-double-critical graph, and assume that $G$ is claw-free. Suppose $|N(x) \cap N(y)| \ge 5$ for all $xy \in E(G)$. 
Then $G \cong K_6$.
\end{lemma}

\pf
We prove this Lemma by way of contradiction. Suppose $G \not \cong K_6$. Then, by Lemma \ref{clique}, $K_5 \not \subset G$. 
By Lemma \ref{neighbor} and Lemma \ref{four}, we may assume that $|N(x) \cap N(y)| = 5$ for all $xy \in E(G)$.
Let $N(x)\cap N(y) = \{v_1,v_2,v_3,v_4,v_5\}$, and let $V_1,V_2,V_3,V_4$ be the color classes of a 4-coloring of $G-x-y$. By Lemma \ref{existence}, we may assume that $v_i \in V_i$, $i \in \{1,2,3,4\}$ and $v_5 \in V_1$.
By Lemma \ref{neighbor}, $\{v_1,v_2,v_3,v_4,v_5\}$ induces a $C_5$ in $G$.
Without loss of generality, assume $v_1v_2,v_1v_4,v_2v_5,v_3v_4,v_3v_5 \in E(G)$ and $v_1v_3,v_2v_3,v_2v_4,v_4v_5 \not\in E(G)$.

Since $|N(x) \cap N(v_5)| = 5$, there exist $a,b \in (N(x) \cap N(v_5)) \setminus \{v_1,v_2,v_3,v_4,v_5,y\}$.
Similarly since $|N(y) \cap N(v_5)| = 5$, there exist $c,d \in (N(y) \cap N(v_5)) \setminus \{v_1,v_2,v_3,v_4,v_5,x\}$.
Note that $a,b,c,d \in (V_2 \cup V_3 \cup V_4) \setminus \{v_2,v_3,v_4\}$.
Since $N(x) \cap N(y) = \{v_1,v_2,v_3,v_4,v_5\}$, $a,b\not\in N(y)$ and $c,d \not\in N(x)$. Hence $a,b,c,d$ are pairwise distinct. Moreover $ab \in E(G)$ to avoid the claw $(x;a,b,y)$ in $G$, and $cd \in E(G)$ to avoid the claw $(y;c,d,x) \in E(G)$.
By Lemma \ref{easy} (applied to $v_5$ and $V_i$, for $i \in \{2,3,4\}$) and by the symmetry between $V_2$ and $V_3$, we may assume that $b, d \in V_4$, $a\in V_3$, and $c\in V_2$.






Since $|N(y)\cap N(v_3)| \ge 5$, there exist $z_1,z_2 \in (N(y) \cap N(v_3)) \setminus \{v_1,v_2,v_3,v_4,v_5,x,y\}$. Note that $z_i \not\in V_1$ for $i \in \{1,2\}$; otherwise $(y;v_1,v_5,z_i)$ induces a claw in $G$. Clearly $z_1,z_2 \not\in V_3$ since $V_3$ is independent. 

We claim that $d \in \{z_1,z_2\}$. Suppose otherwise, $d \not\in \{z_1,z_2\}$. In this case $z_i \not\in V_4$ for $i \in \{1,2\}$ to avoid the claw $(y;d,v_4,z_i)$ in $G$. Therefore $z_1,z_2 \in V_2$, then $(y;v_2,z_1,z_2)$ induces a claw in $G$, a contradiction. 

Since $d \in \{z_1, z_2\}$, we have $dv_3 \in E(G)$.
By a similar argument considering $|N(x) \cap N(v_3)| \ge 5$, $bv_3 \in E(G)$.
Thus $(v_3;b,d,v_4)$ induces a claw in $G$, a contradiction.
\qed



\end{document}